%% file: master.tex
\documentclass[11pt]{article}
\usepackage{times}
\usepackage[hyperindex=true,pageanchor=true,hyperfigures=true,backref=false]{hyperref} 
\usepackage{amsmath}
\usepackage{amssymb}
\usepackage{url}
\usepackage{dsfont} 
\DeclareSymbolFontAlphabet{\Bbb}{AMSb}
\newlength{\myleftmargin}
\input{macros}

\input{local-macros}
\title{A Short Note on the Comparison of Interpolation Widths, Entropy Numbers, and Kolmogorov Widths}

\author{
Ingo Steinwart
}

\begin{document}

\maketitle

\begin{abstract}
We compare the Kolmogorov and entropy numbers of compact operators mapping 
from a Hilbert space into a Banach space. We then apply these general findings 
to embeddings between reproducing kernel Hilbert spaces and $L_\infty(\mu)$. Here 
we provide a sufficient condition for a 
gap of the order $n^{1/2}$ between
the associated  interpolation and Kolmogorov $n$-widths. Finally, we show that in the 
multi-dimensional 
Sobolev case, this gap actually occurs between the Kolmogorov and approximation widths.
\end{abstract}

\input{intro}
\input{prelim}

\input{carls}

\input{main-results}

\input{examples}

\noindent
\textbf{Acknowledgement.}
I deeply thank R.~Schaback, for pointing me to the question regarding the
gap between Kolmogorov and interpolation widths as well as for his feedback and 
his patience at the Shanghai airport.
I also thank G.~Santin for giving valuable hints.

{\small


}

\end{document}

%% file: macros.tex
\usepackage{amsmath}
\usepackage{amssymb}
\usepackage{amsthm}
\usepackage{color}
\DeclareSymbolFontAlphabet{\Bbb}{AMSb}

\newtheorem{theorem}{Theorem}[section]
\newtheorem{lemma}[theorem]{Lemma}

\newtheorem{corollary}[theorem]{Corollary}
\newenvironment{proofof}[1]{\noindent{\bf Proof of #1:}}{\qed\medskip}

\newlength{\fixboxwidth}
\definecolor{darkgreen}{rgb}{0,0.6,0}

\newcommand{\ca}[1]{{\cal #1}}

\newcommand{\N}{\mathbb{N}}

\newcommand{\R}{\mathbb{R}}
\newcommand{\Rd}{\mathbb{R}^d}

\renewcommand{\a}{\alpha}
\renewcommand{\b}{\beta}
\newcommand{\g}{\gamma}

\renewcommand{\d}{\delta}

\newcommand{\e}{\varepsilon}

\newcommand{\lb}{\lambda}

\newcommand{\s}{\sigma}

\renewcommand{\P}{\Phi}

\DeclareMathOperator{\ran}{ran}
\DeclareMathOperator{\rank}{rank}

\DeclareMathOperator{\codim}{codim}

\DeclareMathOperator{\id}{id}

\newcommand{\snorm}[1]{\Vert #1 \Vert}

\newcommand{\bnorm}[1]{\Bigl\Vert \, #1 \, \Bigr\Vert}

\newcommand{\Lx}[2]{{L_{#1}(#2)}}

\newcommand{\tridia}[6]{
\setlength{\unitlength}{1ex}
\begin{picture}(60,20)
\put(10,14){\makebox(9,2)[rt]{$#1$}}
\put(41,14){\makebox(9,2)[lt]{$#2$}}
\put(25.5,2){\makebox(9,2){$#3$}}
\put(20,15){\vector(1,0){20}}
\put(20,14){\vector(1,-1){9}}
\put(31,5){\vector(1,1){9}}
\put(27,16){\makebox(6,2){$#4$}}
\put(17.5,8){\makebox(6,2)[rt]{$#5$}}
\put(37.5,8){\makebox(6,2)[lt]{$#6$}}
\end{picture}}

\newcommand{\quadiasn}[8]{
\setlength{\unitlength}{1ex}
\begin{picture}(60,20)
\put(20,15){\makebox(9,2)[rt]{$#1$}}
\put(45,15){\makebox(9,2)[lt]{$#2$}}
\put(20,2){\makebox(9,2)[rt]{$#3$}}
\put(45,2){\makebox(9,2)[lt]{$#4$}}
\put(30,16){\vector(1,0){14}}
\put(28,14){\vector(0,-1){9}}
\put(46,5){\vector(0,1){9}}
\put(30,3){\vector(1,0){14}}
\put(34,17){\makebox(6,2){$#5$}}
\put(21,8.5){\makebox(6,2)[rt]{$#6$}}
\put(47,8.5){\makebox(6,2)[lt]{$#7$}}
\put(34,0){\makebox(6,2){$#8$}}
\end{picture}}

%% file: intro.tex
\section{Introduction}\label{sec:intro}

Let $(X,\ca A,\mu)$ be a measure space and $H$ be a 
 a reproducing kernel Hilbert space (RKHS) over $X$. Moreover, assume that the kernel $k$ of $H$ is measurable 
and that for all $p\in [2,\infty]$, the map $I_{k,\mu}:H\to L_p(\mu)$ defined by $I_{k,\mu}f:=[f]_\sim$, where 
$[f]_\sim$ denotes the $\mu$-equivalence class of $f$ in $L_p(\mu)$, is compact.
Now consider the
 linear interpolation $n$-width of $H$ in $\Lx 2\mu$, that is 
\begin{displaymath}
   I_n\bigl(H, L_p(\mu)\bigr)
:=
 \inf_{D\subset X, |D| \leq n} \biggl( \int_X  \sup_{f\in B_H} \bigl|f(x) - A_Df(x)   \bigr|^p  \, d\mu(x)  \biggr)^{1/p}\, ,
\end{displaymath}
with the usual modification for $p=\infty$.
Here, $A_D:H\to H$ is the bounded linear operator defined by $A_Df(x) := \sum_{i=1}^n \a_i^*(x) f(x_i)$,
where $D=(x_1,\dots,x_n)$ and $\a^*(x)\in \R^n$ is the unique solution of 
\begin{displaymath}
   \a^*(x) = \arg\min_{\a\in \R^n}  \bnorm{\d_x - \sum_{i=1}^n \a_i\d_{x_i}}_{H'}^2\, .
\end{displaymath}
For later use we note that we always have
\begin{equation}\label{In-lower-bound}
 \inf_{D\subset X, |D| \leq n} \sup_{f\in B_H} \snorm{f-A_Df}_{L_p(\mu)} \leq  I_n\bigl(H, L_p(\mu)\bigr)
\end{equation}
and equality holds in the extreme case $p=\infty$.
Moreover, consider the classical Kolmogorov $n$-width 
\begin{displaymath}
 d_n\bigl(H,\Lx p\mu\bigr) = \inf_{F_{n}\subset \Lx p\mu}\sup_{f\in B_H}\inf_{g\in F_n}\snorm{f-g}_{\Lx p\mu}\, ,
\end{displaymath}
where the left most infimum runs over all subspaces $F_n$ of $\Lx p\mu$ with $\dim F_n\leq n$. 
Note that the lower bound of  $I_n$ in \eqref{In-lower-bound}
measures, how well $f$ can be approximated by a very particular  linear and $n$-dimensional 
scheme, whereas $d_n$ measures how well $f$ can be approximated by the best $n$-dimensional scheme.
Consequently, the approximation $n$-width 
\begin{displaymath}
a_n\bigl(H, L_p(\mu)\bigr) := \inf_{A:H\to \Lx p \mu} \sup_{f\in B_H} \snorm{f-Af}_{L_p(\mu)}\, ,
\end{displaymath}
where the infimum is taken over all bounded linear operators $A:H\to \Lx p \mu$ with $\rank A\leq n$, satisfies
$d_n(H, L_p(\mu)) \leq a_n(H, L_p(\mu)) \leq I_n(H, L_p(\mu))$. 

In the Hilbert space case, that is, $p=2$, these quantities are well understood. Indeed, the general 
theory of $s$-numbers \cite{Pietsch87} shows, see e.g.~Section \ref{sec:prelim}, that 
\begin{equation}\label{equal-hilbert}
d_n\bigl(H, L_2(\mu)\bigr) = a_n\bigl(H, L_2(\mu)\bigr) = \sqrt{\lb_{n+1}}\, ,
\end{equation} 
where $(\lb_n)$ denotes the (extended and) ordered sequence of eigenvalues of the integral operator 
$T_k: L_2(\mu) \to L_2(\mu)$ associated with the kernel $k$. Moreover, if $H$ is a Sobolev space, then 
$I_n(H, L_p(\mu))$ shares the asymptotic behavior of \eqref{equal-hilbert} and this can actually be achieved by 
taking quasi-uniform points
$D\subset X$, see 
\cite{ScWe02a}. Unfortunately, the situation changes in the other extreme, namely $p=\infty$.
Indeed, if $\mu$ is a finite measure, then \eqref{equal-hilbert} immediately yields
\begin{displaymath}
    \sqrt{\lb_{n+1}} \leq \sqrt{\mu(X)} \, d_n\bigl(H, L_\infty(\mu)\bigr)\, ,
\end{displaymath}
while \cite[Theorem 3]{SaSc16a} shows that 
\begin{equation}\label{in-lower-bound}
    \sqrt{\sum_{i=n+1}^\infty\lb_{i}} \leq \sqrt{\mu(X)} \, I_n\bigl(H, L_\infty(\mu)\bigr)\, ,
\end{equation}
and in the Sobolev case, this lower bound is matched 
by an upper bound of the same asymptotic behavior, see \cite{ScWe02b}.
In the case of an algebraic decay of the eigenvalues, it is not hard to see that there is a gap of the order $n^{-1/2}$
between the \emph{lower bounds} for $d_n(H, L_\infty(\mu))$ and $I_n(H, L_\infty(\mu))$, and this 
naturally raises the question, whether this gap actually occurs between the quantities of interest, that is, between
$d_n(H, L_\infty(\mu))$ and $I_n(H, L_\infty(\mu))$. So far, a positive answer only exists for the 
1-dimensional Sobolev case, see \cite{Pinkus85}. The goal of this note is to provide a positive answer 
in a more general framework. To be more precise, we  show that for algebraically decaying eigenvalues we have 
$d_n(H, L_\infty(\mu)) \asymp \sqrt{\lb_{n+1}}$ if and only if the entropy numbers of the embedding 
$I_{k,\mu}:H\to \Lx \infty\mu$ behave like $\sqrt{\lb_{n+1}}$. Using 
\eqref{in-lower-bound} this characterization gives a sufficient condition for the existence of the gap.
In addition, we present a result that highlights the role of the eigenfunctions of $T_k$.
For the multi-dimensional Sobolev case 
we then show with the help of  well-known asymptotics of the entropy and approximation
numbers
that the gap $n^{-1/2}$ actually occurs between $d_n(H, L_\infty(\mu))$ and $a_n(H, L_\infty(\mu))$, that is, 
between arbitrary $n$-dimensional approximation and  linear $n$-dimensional approximation.
In addition, the cases $p\in (2,\infty)$ are treated simultaneously.

The rest of this note is organized as follows: In Section \ref{sec:prelim} we recall the definition of 
entropy numbers and also introduce some examples of $s$-number scales. Section \ref{sec:carls}
summarizes the relationship between entropy numbers and the different $s$-number scales. 
In Section \ref{sec:main} 
two general results comparing entropy and Kolmogorov numbers of compact operators are presented
and based upon these results the RKHS situation is investigated in more detail.
In Section \ref{sec:examples} we then apply these findings to the multi-dimensional Sobolev case. 

%% file: prelim.tex
\section{Preliminaries: Entropy Numbers, $\mathbf{\mathit{s}}$-Numbers, and Eigenvalues}\label{sec:prelim}

We write $a_n \prec b_n$ for two positive sequences $(a_n)$ and $(b_n)$
if there exists a constant $c\in (0,\infty)$ such that $a_n \leq c b_n$ for all $n\geq 1$.
Similarly, we write $a_n \asymp b_n$ if both $a_n \prec b_n$ and $b_n\prec a_n$.
Finally, a positive sequence is called regular if there exists a constant $c\in (0,\infty)$ such that
 $a_n \leq c a_{2n}$ and 
$a_m\leq ca_n$ for all $1\leq m\leq n$. Probably the most interesting examples of regular 
sequences are $a_n = n^{-p} (1 + \ln n)^{-q}$ for $p> 0$ and $q\in \R$, or $p=0$ and $q>0$.


Given a Banach space $E$, we denote its closed unit ball by $B_E$ and its dual by $E'$.
Moreover, we write  $I_F:F\to \ell_\infty(B_{F'})$ for  the canonical embedding
and 
$Q_E:\ell_1(B_E) \to E$ for the canonical surjection. Furthermore, we write $E\hookrightarrow F$
if $E\subset F$ and the inclusion map is continuous.
Finally, the adjoint of a bounded linear 
operator $S$  acting between two Hilbert spaces is denoted by $S^*$.

Now, let $E$ and $F$ be Banach spaces and $T:E\to F$ be a bounded, linear operator.
Then the $n$-th \emph{(dyadic) entropy number} of $T$ is defined by 
\begin{displaymath}
 e_n(T) := \inf\Bigl\{\e>0: \exists x_1,\dots,x_{2^{n-1}}\in F: TB_E \subset \bigcup_{i=1}^{2^{n-1}} x_i + \e B_F\Bigr\}\, .
\end{displaymath}
Some elementary properties of entropy numbers can be found in \cite[Chapter 1]{CaSt90}. In particular, we have 
$e_n(T)\to 0$ if and only if $T$ is compact. Since $T$ is compact if and only if its dual $T'$ is compact, this 
immediately raises 
 the question how the entropy numbers of $T$ and $T'$ are related to each other. This question, known
as the duality problem for entropy numbers has, so far, no complete answer. Partial answers, however, 
do exist. The one we will need  is the following inequality taken from \cite{BoPaSzTJ89a}
\begin{equation}\label{dual-entr}
\frac 1 {d_p}\ \sup_{k\leq n} k^{1/p}\ e_k(T) \leq 
\sup_{k\leq n} k^{1/p}\ e_k(T') \leq d_p\ \sup_{k\leq n} k^{1/p}\ e_k(T)\, ,
\end{equation}
which holds for all $n\geq 1$ and all compact $T:E\to F$,
whenever $E$ or $F$ are $B$-convex. Here, $d_p\in (0,\infty)$ is a constant, which depends on $p\in (0,\infty)$ 
and the geometry of the involved spaces $E$ and $F$, but which is 
independent of both $n$ and $T$. 
Moreover, recall from e.g.~\cite[Theorem 13.10]{DiJaTo95} that a Banach space is $B$-convex if and only if it 
has non-trivial type. In particular, Hilbert spaces are $B$-convex, and so are the spaces 
$L_p(\mu)$ for $p\in (1,\infty)$ since these spaces have type $\min\{2,p\}$, see e.g.~\cite[Chapter 11]{DiJaTo95}.
Moreover, if $E$ or $F$ is a Hilbert space, it was shown in \cite{TJ87a}
that we may choose  $d_p=32$ for all $p\in (0,\infty)$.
Finally note that 
from the inequalities in \eqref{dual-entr} 
we can derive the following equivalences, which hold for all regular sequences  $(\a_n)$
and all compact operators $T$:
\begin{align*}
   e_n(T) \prec \a_n \qquad &\Longleftrightarrow   \qquad e_n(T') \prec \a_n \\
 e_n(T) \asymp \a_n \qquad &\Longleftrightarrow   \qquad e_n(T') \asymp \a_n \, .
\end{align*}
For a proof, 
which is based on a little trick originating from  Carl
\cite{Carl85a},
 we refer to the proof of 
 \cite[Corollary 1.19]{Steinwart00b} or, in a slightly simplified version, 
to the proof of  \cite[Proposition 2]{Steinwart00a}.

Besides entropy numbers, we are also interested in some so-called $s$-numbers. Namely, if $T:E\to F$ is a bounded linear 
operator, we are interested in the 
 $n$-th {\em approximation number} of $T$, defined by
\begin{displaymath}
a_n(T) := \inf \{\ \snorm{T-A}\ | \ A:E\to F \mbox{ bounded, linear with }
\rank A < n\ \}\, ,
\end{displaymath}
in the $n$-th {\em Gelfand number} of $T$ defined by
\begin{displaymath}
c_n(T)  :=  \inf \{\ \snorm{T I_{E_0}^E} \ : \ 
\mbox{$E_0$ subspace of $E$ with }\codim {E_0} <n\ \}\, ,
\end{displaymath}
where $I_{E_0}^E$ denotes the canonical inclusion of $E_0$ into $E$, and in
the $n$-th {\em Kolmogorov number} of $T$ defined by
\begin{displaymath}
d_n(T)  :=  \inf \{\ \snorm{Q_{F_0}^F T} \ : \
\mbox{$F_0$ subspace of $F$ with } \dim F_0 <n\ \}\, ,
\end{displaymath}
where $Q_{F_0}^F$ denotes the canonical surjection from the Banach space $F$ onto
the quotient space $F/F_0$. 
Recall from \cite[Proposition 2.2.2]{CaSt90} that the 
 latter quantity can also be expressed by 
\begin{displaymath}
  d_n(T) = \inf\Bigl\{ \e>0: \exists F_\e \mbox{ subspace of } F \mbox{ with } \dim F_\e < n \mbox{ and } TB_E \subset  F_\e+\e B_F         \Bigr\}\, ,
\end{displaymath}
and consequently, we have
\begin{displaymath}
 d_{n+1}(T) = \inf_{F_{n}\subset F}\sup_{y\in TB_E}\inf_{z\in F_n}\snorm{y-z}_F\, ,
\end{displaymath}
where the left most infimum runs over all subspaces $F_n$ of $F$ with $\dim F_n\leq n$. 
In other words, $d_{n+1}(T)$ equals the classical Kolmogorov $n$-width of the set $TB_E$ in $F$, 
cf.~\cite[Chapter 13]{LoGoMa96}, and therefore we have
$d_{n+1}(I_{k,\mu}:H\to \Lx p \mu) = d_n(H,\Lx p \mu)$, where $H$ and $\mu$ are as in the introduction.
In addition, it is not hard to see that we also have $a_{n+1}(I_{k,\mu}:H\to \Lx p \mu) = a_n(H,\Lx p \mu)$,
and consequently we will consider the operator versions in the remaining parts of this note.
Furthermore, 
recall e.g.~from \cite[Theorems 2.3.1 and 2.2.1, and Proposition 2.5.5]{CaSt90} that we always have 
\begin{align*}
 c_n(T) &= a_n(I_FT)\\
 d_n(T) &= a_n(TQ_E) \\ 
 d_n(T')&= c_n(T)\, , 
\end{align*}
and for compact $T$ its dual operator $T'$  additionally
satisfies $c_n(T')= d_n(T)$, see \cite[Proposition 2.5.6]{CaSt90}.
Moreover, the approximation, Gelfand, and Kolmogorov numbers are $s$-numbers in the sense of \cite[Definition 2.2.1]{Pietsch87}, and the same is 
true for the  {\em Tikhomirov numbers} of $T$, which are  defined by 
\begin{displaymath}
t_n(T)  :=  a_n(I_FTQ_E) \, , \qquad \qquad n\geq 1.
\end{displaymath}
In addition, we always have 
\begin{align*}
 t_n(T) & \leq c_n(T)\leq a_n(T) \leq  \sqrt {2n}\,  c_n(T) \\
t_n(T)  & \leq d_n(T)\leq a_n(T) \leq  \sqrt {2n}\, d_n(T)\, , 
\end{align*}
where we note that in both cases the first two inequalities follow from $s$-number properties and 
the right most inequalities can be found in \cite[Propositions 2.4.3 and 2.4.6]{CaSt90}.
In addition, the factor $\sqrt {2n}$ can be sharpened to $1+\sqrt{n-1}$.

The two chains of inequalities above show that the gap between the asymptotic behavior 
of $(a_n)$, $(c_n)$ and $(d_n)$ is at most of the order $\sqrt n$. It is well-known that 
this gap is sometimes attained, see e.g.~Section \ref{sec:examples}, while in other cases
the gap vanishes. 
For example, we have 
\begin{equation}\label{no-gap-an-cn}
 a_n(T) =  c_n(T)\, ,
\end{equation}
if $E$ is a Hilbert space, see \cite[Proposition 2.4.1]{CaSt90}, or $F$ 
has the metric extension property, see \cite[Proposition 2.3.3]{CaSt90},
and 
\begin{equation}\label{no-gap-an-dn}
 a_n(T) =  d_n(T)\, , 
\end{equation}
if $F$ is a Hilbert space, see \cite[Proposition 2.4.4]{CaSt90}, or $E$ has the metric lifting property,
see \cite[Proposition 2.2.3]{CaSt90}. In this respect recall that the spaces $\ell_\infty(J)$ and
$L_\infty(\mu)$, where $\mu$ is  some finite measure, have the metric extension property,
see \cite[p.~60]{CaSt90} and \cite[Theorem 4.14]{DiJaTo95}, respectively.
Moreover, the spaces $\ell_1(I)$ have the metric lifting property, see \cite[p.~51]{CaSt90}.
By combining all these relations we further see that 
we have $t_n(T) = a_n(T)$ 
if either  $E=\ell_1(I)$ and $F=\ell_1(J)$, or $E$ and $F$ are Hilbert spaces. The latter case also follows from a general 
result showing that there is only one $s$-scale for operators between Hilbert spaces, 
see e.g.~\cite[Theorem 2.11.9]{Pietsch87}.

Our next goal is to relate the $s$-numbers introduced above to eigenvalues. To this end, let 
$S:H_1\to H_2$ be a compact operator acting between two Hilbert spaces. Then $S^*S:H_1\to H_1$
is compact, self-adjoint and positive, and therefore the classical spectral theorem shows that 
there is an at most countable family $(\lb_i(S^*S))_{i\in I}$ of eigenvalues of $S^*S$, which in addition
are non-negative and have at most one limit point, namely $0$. In the following, we always assume 
that either $I=\{1,\dots,n\}$ or $I=\N$, and that the eigenvalues are ordered 
decreasingly without excluding (geometric) multiplicities. Then, the \emph{singular numbers} of $S$ are defined by  
\begin{displaymath}
s_i (S):= 
\begin{cases}
 \sqrt{\lb_i(S^*S)} = \lb_i(\sqrt{S^*S}) & \mbox{ if } i\in I\\
0 & \mbox{ if } i\in \N\setminus I \, .
\end{cases}
\end{displaymath}
Recall that this gives $s_i(S) = s_i(S^*)$ for all $i \geq 1$, 
and $s_i(T)=\lb_i(T)$ for all $i\in I$ if $T:H\to H$ is compact, self-adjoint and positive.
Moreover, we have, see e.g.~\cite[Chapter 2.11]{Pietsch87}
\begin{displaymath}
   s_n(S) = a_n(S)
\end{displaymath}
for all $n\geq 1$ and all compact operators $S:H_1\to H_2$ between Hilbert spaces $H_1$, $H_2$.

%% file: carls.tex
\section{Carl's Inequality and some Inverse Versions}\label{sec:carls}

In this section we   recall some inequalities between $s$-numbers and entropy numbers. We begin with 
Carl's inequality, which states that for all $p\in (0,\infty)$, there exists a constant $C_p\in (0,\infty)$
such that for all bounded, linear $T:E\to F$ and all $n\geq 1$, we have 
\begin{equation}\label{carls-inequality}
   \sup_{k\leq n} k^{1/p} e_k(T) \leq C_p \sup_{k\leq n} k^{1/p} a_k(T) \, .
\end{equation}
We refer to \cite[Theorem 3.1.1]{CaSt90}, where it is also shown that a possible value for the constant is
$C_P = 128 (32+16/p)^{1/p}$. Recall from e.g.~\cite[Chapter 1.3]{CaSt90} that entropy numbers are surjective and 
weakly injective, and therefore we have 
\begin{displaymath}
   e_k(T) \leq 2e_k(I_F T Q_E) \leq 2e_k(T) 
\end{displaymath}
for all bounded, linear  $T:E\to F$ and all $k\geq 1$. In particular, we may replace the approximation numbers in 
\eqref{carls-inequality} by the Gelfand, Kolmogorov, or Tikhomirov numbers for the price of an additional 
factor of $2$ in the constant. 
Moreover, 
like for the entropy numbers of $T$ and $T'$,  we further have 
\begin{displaymath}
   a_n(T) \prec \a_n \qquad \Longrightarrow   \qquad e_n(T) \prec \a_n
\end{displaymath}
for all regular sequences $(\a_n)$ and all bounded linear
$T:E\to F$. It is needless to say that the approximation numbers in this implication 
 may be replaced by 
the Gelfand, Kolmogorov, or Tikhomirov numbers.

Let us now recall some inequalities that describe how certain $s$-numbers are dominated by 
entropy numbers. We begin with compact operators $S:H_1\to H_2$ acting between two Hilbert spaces.
Then \cite[Inequality (3.0.9)]{CaSt90} shows 
\begin{equation}\label{an-vs-en}
   a_n(S) \leq 2 e_n(S)
\end{equation}
for all $n\geq 1$. 
By an adaptation of the proof of 
 \cite[Corollary 1.19]{Steinwart00b} we can then see that  \eqref{an-vs-en} in combination with 
\eqref{carls-inequality} leads to
%
the following equivalences, which hold for all regular sequences  $(\a_n)$
and all compact operators $S:H_1\to H_2$ acting between two Hilbert spaces:
\begin{align}\label{Hilbert-aquiv1}
   a_n(S) \prec \a_n \qquad &\Longleftrightarrow   \qquad e_n(S) \prec \a_n \\ \label{Hilbert-aquiv2}
 a_n(S) \asymp \a_n \qquad &\Longleftrightarrow   \qquad e_n(S) \asymp \a_n \, .
\end{align}
Again, the approximation numbers in these equivalences
 may be replaced by 
the Gelfand, Kolmogorov, or Tikhomirov numbers.
Finally, let us consider the compact, self-adjoint and positive 
operator $T:H_1\to H_1$ defined by $S^*S$. Then we have 
\begin{equation}\label{s-squared}
   s_i(T) = \lb_i(T) = \lb_i(S^*S) = s_i^2(S^*)
\end{equation}
if $i\in I$ and $s_i(T) = 0 = s_i^2(S^*)$ if $i\in \N\setminus I$.
The two equivalences above then lead to 
\begin{align}\label{ln-vs-en-1}
   s_n(T) \prec \a_n \qquad &\Longleftrightarrow   \qquad e_n(S^*) \prec \sqrt{\a_n} \\  \label{ln-vs-en-2}
 s_n(T) \asymp \a_n \qquad &\Longleftrightarrow   \qquad e_n(S^*) \asymp \sqrt{\a_n} 
\end{align}
for all regular sequences $(\a_n)$. Note that $s_n$ can be replaced by 
any $s$-number scale, and in particular by the approximation, Gelfand, Kolmogorov,
and Tikhomirov numbers. Moreover, we may replace $e_n(S^*)$ by $e_n(S)$ using the 
duality results for entropy numbers mentioned above. 

Let us now consider the situation in which only one of the involved spaces  is a Hilbert space, that is, 
we consider compact operators of the form $S:E\to H$ or 
$S:H\to F$, where $H$ is a Hilbert space and $E$ or $F$ is an arbitrary Banach space.
Then \eqref{carls-inequality} still holds, but in general, we may no longer have \eqref{an-vs-en}.
To compare the $s$-numbers of $T$ to the entropy numbers of $T$, we thus need a surrogate for  
 \eqref{an-vs-en}. Fortunately, there are a few such results. 
For example, \cite[Lemme 1]{PaTJ85a} shows that there exist constants $A,B\in (0,\infty)$ such that 
for all compact $S:E\to H$ and all $n\geq 1$ we have 
\begin{equation}\label{PT-1}
   n^{1/2} c_n(S) \leq B\sum_{k> A n}k^{-1/2} e_k(S)\, .
\end{equation}
With the help of this inequality it is easy to show that for all $p\in (0,2)$ there exists 
another constant $B_p\in (0,\infty)$ such that 
\begin{equation}\label{PT-2}
    n^{1/p} c_n(S) \leq B_p \sup_{k> An} k^{1/p} e_k(S)
\end{equation}
for all compact $S:E\to H$ and all $n\geq 1$. We refer to the very short proof of 
\cite[Th\'eor\`eme A]{PaTJ85a}. Complementary, \cite[Theorem 5.12]{Steinwart00b} shows that 
for all $p\in (2,\infty)$ there exists a constant $K_p\in (0,\infty)$ such that 
\begin{equation}\label{carls-inequality-inverse}
  \sup_{k\leq n} k^{1/p} t_k(S) \leq K_p \,\sup_{k\leq n} k^{1/p} e_k(S) \, .
\end{equation}
for all compact operators $S:E\to H$ or $S:H\to F$ and all $n\geq 1$.
Last but not least we like to mention that \cite[Theorem 6]{Carl85a} 
showed an inequality of the form \eqref{carls-inequality-inverse}
with $t_k$ replaced by $d_k$ or $c_k$ for all $p\in (0,\infty)$ and all compact  $S:E\to F$
for which $E$ and $F'$ are type 2 spaces.

%% file: main-results.tex
\section{Main Results}\label{sec:main}

The goal of this section is to compare the entropy and Kolmogorov numbers of the embedding
$I_{k,\mu}:H\to \Lx\infty\mu$. 
To this end, our first auxiliary  result
combines Carl's inequality with its inversed versions mentioned in Section \ref{sec:carls}.

\begin{lemma}\label{new-equiv}
Let $H$ be a Hilbert space, $F$ be a Banach space 
$S:H\to F$ be a compact operator, and $p\in (0,2)$. Then, the following equivalence holds: 
\begin{align}\label{new-equiv-h1}
      d_n(S) \prec  n^{-1/p}  \qquad &\Longleftrightarrow   \qquad  e_n(S) \prec  n^{-1/p} \, .
\end{align}
Moreover, if $F$ has the metric extension property, the equivalence is also true for $p\in (2,\infty)$, 
and in addition, we have
\begin{align}\label{new-equiv-h2}
      d_n(S) \asymp  n^{-1/p}  \qquad &\Longleftrightarrow   \qquad  e_n(S) \asymp  n^{-1/p} \, .
\end{align}
Finally, if $F'$ has type 2, then \eqref{new-equiv-h1} and \eqref{new-equiv-h2} hold for all $p\in (0,\infty)$. 
\end{lemma}

\begin{proofof}{Lemma \ref{new-equiv}}
Independent of $p$ and $F$, 
the implication ``$\Rightarrow$'' in \eqref{new-equiv-h1}
 is a direct consequence of Carl's inequality \eqref{carls-inequality}.
For the proof of the converse implication 
we first consider the case $p\in (0,2)$.
By \eqref{dual-entr} we then know that $e_n(S') \prec  n^{-1/p}$, and consequently \eqref{PT-2}
shows that $c_n(S') \prec  n^{-1/p}$. Using $c_n(S') = d_n(S)$, which holds for compact operators $S$, 
we then obtain the assertion.
In the case $p\in (2,\infty)$, we conclude by \eqref{carls-inequality-inverse} that 
$t_n(S) \prec  n^{-1/p}$. Moreover, $F$ has the metric extension property, and therefore 
we have $c_n(SQ_E) = a_n(SQ_E)$ by \eqref{no-gap-an-cn}.
This leads to
\begin{displaymath}
   t_n(S) = a_n(I_FSQ_E)=c_n(SQ_E) = a_n(SQ_E) = d_n(S)\, ,
\end{displaymath}
and hence we find $d_n(S) \prec  n^{-1/p}$. In addition, \eqref{new-equiv-h2} follows from 
combining  \eqref{carls-inequality} and  \eqref{carls-inequality-inverse} as in 
the proof of  \cite[Proposition 2]{Steinwart00a}. Finally, the last assertion can be shown analogously using 
 \cite[Theorem 6]{Carl85a} instead of \eqref{carls-inequality-inverse}.
\end{proofof}

Note that the equivalences obtained in 
Lemma \ref{new-equiv} also holds for regular sequences of the form $\a_n = n^{-1/p} (\log n)^\b$, where 
$p$ satisfies the constraints of Lemma \ref{new-equiv} and $\b\in \R$. Indeed, for the second and third 
case this 
can be deduced from \eqref{carls-inequality} and \eqref{carls-inequality-inverse}, respectively 
\cite[Theorem 6]{Carl85a}, while in the first  case 
this follows from  \eqref{carls-inequality} , \eqref{PT-2}, and \cite[G.3.2]{Pietsch87}.

Clearly,
Lemma \ref{new-equiv} in particular holds for compact operators $S:H\to L_p(\mu)$. Our next result shows that 
for some spaces $L_p(\mu)$ even more information can be obtained.

\begin{theorem}\label{main1}
Let $H$ be a Hilbert space, $\mu$ be a finite measure, and $p\in [2,\infty]$. Assume that we have 
a compact operator $S:H\to L_p(\mu)$ such that 
\begin{equation}\label{hilbert-condition}
   e_n\bigl( S:H\to L_2(\mu) \bigr) 
\asymp n^{-1/\a}
\end{equation}
for some $\a\in (0,2)$. Then, for all $q\in [2,p]$, the following equivalence hold: 
\begin{align*}
d_n\bigl( S:H\to L_q(\mu) \bigr) \asymp  n^{-1/\a}  \qquad &\Longleftrightarrow   \qquad  e_n\bigl( S:H\to L_q(\mu) \bigr) \asymp  n^{-1/\a} \, .
\end{align*}
\end{theorem}

\begin{proofof}{Theorem \ref{main1}}
``$\Rightarrow$'': By Lemma \ref{new-equiv}, or more precisely,
Carl's inequality,  we already know  that $e_n( S:H\to L_q(\mu) )\prec n^{-1/\a}$.
Moreover, using $L_q(\mu) \hookrightarrow L_2(\mu)$ we find
\begin{displaymath}
  n^{-1/\a} \asymp  e_n\bigl( S:H\to L_2(\mu) \bigr) 
\leq \snorm{\id: L_q(\mu) \to L_2(\mu)} \, e_n\bigl( S:H\to L_q(\mu) \bigr)\, ,
\end{displaymath}
and thus $e_n( S:H\to L_q(\mu) )\asymp n^{-1/\a}$.

``$\Leftarrow$'': By Lemma  \ref{new-equiv},  we already know  that $d_n( S:H\to L_q(\mu) )\prec n^{-1/\a}$.
Moreover, by \eqref{hilbert-condition}, \eqref{Hilbert-aquiv2}, 
and \eqref{no-gap-an-dn}
we obtain  
$d_n( S:H\to L_2(\mu) )\asymp n^{-1/\a}$, and hence we find  
\begin{displaymath}
  n^{-1/\a} \asymp  d_n\bigl( S:H\to L_2(\mu) \bigr) 
\leq \snorm{\id: L_q(\mu) \to L_2(\mu)} \, d_n\bigl( S:H\to L_q(\mu) \bigr)\, ,
\end{displaymath}
that is $d_n( S:H\to L_q(\mu) )\asymp n^{-1/\a}$.
\end{proofof}

Note that the entropy numbers in condition \eqref{hilbert-condition}
can be replaced by the Kolmogorov numbers. Indeed, \eqref{Hilbert-aquiv2} shows that 
\eqref{hilbert-condition} is equivalent to $a_n( S:H\to L_2(\mu))\asymp n^{-1/\a}$, and 
since we further have $a_n(S) = d_n(S)$, we see that 
 condition \eqref{hilbert-condition}
can be replaced by 
\begin{equation}\label{hilbert-condition-alter}
   d_n\bigl( S:H\to L_2(\mu) \bigr)\asymp n^{-1/\a}\, .
\end{equation}
In addition, if $H$ is an RKHS with kernel $k$ and 
$T_k$ denotes the integral operator associated with $k$,
then \eqref{hilbert-condition-alter}, or 
\eqref{hilbert-condition}, can  be replaced by 
\begin{equation}\label{hilbert-condition-EV}
   \lb_n\bigl( T_k:L_2(\mu)\to L_2(\mu) \bigr)\asymp n^{-2/\a}
\end{equation}
with the help of \eqref{s-squared}. The following corollary summarizes our findings in this 
situation in view of the gap discussed in the introduction.

\begin{corollary}\label{gap-1}
Let $H$ be  an RKHS of a bounded measurable kernel $k$ on $(X,\ca A)$ and $\mu$ be a finite measure
on the $\s$-algebra $\ca A$. If, in addition, we have
\begin{displaymath}
   e_n\bigl(I_{k,\mu}:H\to \Lx 2 \mu\bigr) \asymp  e_n\bigl(I_{k,\mu}:H\to \Lx \infty \mu\bigr)\asymp n^{-1/\a}
\end{displaymath}
for some $\a\in (0,2)$, then we have $d_n(H, \Lx \infty \mu)\asymp n^{-1/\a}$ and $n^{-1/\a + 1/2} \prec  I_n(H, L_\infty(\mu))$.
\end{corollary}

\begin{proofof}{Corollary \ref{gap-1}}
The behavior $d_n(I_{k,\mu}:H\to \Lx \infty \mu)\asymp n^{-1/\a}$ follows from Theorem \ref{main1}. Moreover, 
we know $\lb_i \asymp  i^{-2/\a}$ by  \eqref{s-squared} and \eqref{ln-vs-en-2}, and therefore, \eqref{in-lower-bound}
shows
\begin{displaymath}
  n^{-1/\a + 1/2} \prec
  \sqrt{\sum_{i=n+1}^\infty i^{-2/\a}} \prec   \sqrt{\sum_{i=n+1}^\infty\lb_{i}} \leq \sqrt{\mu(X)} \, I_n\bigl(H, L_\infty(\mu)\bigr)\, ,
\end{displaymath}
that is, we have shown the second assertion, too.
\end{proofof}

Our last result in this section shows that in the RKHS case 
and $q=\infty$ the asymptotic behavior
 $ e_n(I_{k,\mu}:H\to \Lx 2 \mu) \asymp  n^{-1/\a}$ is inherited from certain 
interpolation spaces between $H$ and $\Lx 2 \mu$. For its formulation we need
the scale of interpolation spaces of the real method, see e.g.~\cite[Chapter 5]{BeSh88}, as well as 
the notation $[H]_\sim := \{[f]_\sim: f\in H\}$.

\begin{theorem}\label{interpol}
Let $H$ be  an RKHS of a bounded measurable kernel $k$ on $(X,\ca A)$ and $\mu$ be a finite measure
on the $\s$-algebra $\ca A$ such that $\ca A$ is $\mu$-complete and that 
\begin{equation}\label{entrop-assump}
   e_n\bigl(I_{k,\mu}:H\to \Lx 2 \mu\bigr) \asymp n^{-1/\a}
\end{equation}
for some $\a\in (0,2)$. 
 In addition assume that
the interpolation space $[\Lx 2 \mu,  [H]_\sim]_{\b,2}$ is compactly embedded into $\Lx \infty \mu$
 for some $\b\in (\a/2,1]$ with 
\begin{equation}\label{entrop-assump-2}
   e_n\bigl([\Lx 2 \mu, [H]_\sim]_{\b,2}\hookrightarrow \Lx \infty\mu\bigr) \prec n^{-\b/\a}
\end{equation}
Then we   have 
\begin{equation}\label{interpol-result}
   e_n\bigl([\Lx 2 \mu, [H]_\sim]_{\g,2}\hookrightarrow \Lx \infty\mu\bigr) \asymp  e_n\bigl([\Lx 2 \mu, [H]_\sim]_{\g,2}\hookrightarrow \Lx 2\mu\bigr)  \asymp n^{-\g/\a}
\end{equation}
for all $\g\in [\b,1]$, as well as $e_n(I_{k,\mu}:H\to \Lx \infty \mu) \asymp n^{-1/\a}$.
\end{theorem}

Before we prove this theorem we note that the interpolation spaces 
$[\Lx 2 \mu, [H]_\sim]_{\g,2}$ can be identified as RKHSs provided that the 
 assumptions of Theorem \ref{interpol} are satisfied. For details we
refer to the  last part of the following proof.

\vspace*{1ex}
\begin{proofof}{Theorem \ref{interpol}}
We first consider the case $\b\in (\a/2,1)$.
  Since $\g\geq \b$, we then know $[\Lx 2 \mu, [H]_\sim]_{\g,2}\hookrightarrow [\Lx 2 \mu, [H]_\sim]_{\b,2}$
by e.g.~\cite[Theorems 4.3 and 4.6]{StSc12a}, and consequently
we have the following diagram of bounded linear embeddings:

\tridia {[\Lx 2 \mu, [H]_\sim]_{\g,2}}{\Lx \infty \mu}{[\Lx 2 \mu, [H]_\sim]_{\b,2}} {}{}{}

\noindent
The multiplicativity  of entropy numbers thus yields
\begin{align}\nonumber
& 
   e_{2n}\bigl([\Lx 2 \mu, [H]_\sim]_{\g,2}\hookrightarrow \Lx \infty\mu\bigr)\\ \label{interpol-h1}
&\leq 
 e_n\bigl([\Lx 2 \mu, [H]_\sim]_{\g,2}\hookrightarrow [\Lx 2 \mu, [H]_\sim]_{\b,2}\bigr)
\cdot
 e_n\bigl([\Lx 2 \mu, [H]_\sim]_{\b,2}\hookrightarrow \Lx \infty\mu\bigr)\, .
\end{align}
Now recall from \cite[Equation (36) and Theorem 4.6]{StSc12a} that 
\begin{equation}\label{interpol-represent}
[\Lx 2 \mu, [H]_\sim]_{\g,2} = \biggl\{ \sum_{i\in I}  a_i \lb_i^{\g/2} [e_i]_\sim: (a_i)\in \ell_2(I)    \biggr\}\, ,
\end{equation}
where $(\lb_i)$  is the sequence of eigenvalues  of the integral operator 
$T_k:\Lx 2 \mu\to \Lx 2 \mu$ and $([e_i]_\sim)$ is a corresponding ONS of eigenfunctions.
Moreover, the system $(\lb_i^{\g/2} [e_i]_\sim)$ is an ONB of $[\Lx 2 \mu, [H]_\sim]_{\g,2}$
with respect to an equivalent Hilbert space norm on $[\Lx 2 \mu, [H]_\sim]_{\g,2}$.
Consequently, we have the following diagram for the embedding $[\Lx 2 \mu, [H]_\sim]_{\g,2}\hookrightarrow [\Lx 2 \mu, [H]_\sim]_{\b,2}$:

\quadiasn {[\Lx 2 \mu, [H]_\sim]_{\g,2}}{[\Lx 2 \mu, [H]_\sim]_{\b,2}} {\ell_2(I)}{\ell_2(I)}  {} {\P_\g} {\P_\b^{-1}}  {D_\Lambda^{(\g-\b)/2}}

\vspace*{1ex}
\noindent 
where $\P_\g$ and $\P_\b$ are the coordinate mappings and $D_\Lambda^{(\g-\b)/2}$ is the diagonal operator
associated to the sequence $(\lb_i^{(\g-\b)/2})$. By \eqref{entrop-assump}, \eqref{s-squared}, and \eqref{ln-vs-en-2}
we conclude that $\lb_i^{(\g-\b)/2} \asymp i^{-(\g-\b)/\a}$.
Using \cite[Proposition 2]{Carl81a}, which estimates entropy numbers of diagonal operators, and the diagram above,
we thus find  
\begin{displaymath}
    e_n\bigl([\Lx 2 \mu, [H]_\sim]_{\g,2}\hookrightarrow [\Lx 2 \mu, [H]_\sim]_{\b,2}\bigr) \prec n^{-(\g-\b)/\a}\, .
\end{displaymath}
Combining this with \eqref{interpol-h1}, \eqref{entrop-assump-2}, and the fact that $\mu$ is a finite measure we obtain
\begin{displaymath}
 e_n\bigl([\Lx 2 \mu, [H]_\sim]_{\g,2}\hookrightarrow \Lx 2\mu\bigr) \prec  e_n\bigl([\Lx 2 \mu, [H]_\sim]_{\g,2}\hookrightarrow \Lx \infty\mu\bigr) \prec  n^{-\g/\a}\, .
\end{displaymath}
To establish the lower bound, we recall from \cite[Proposition 4.2 and Theorems 5.3 and 4.6]{StSc12a} that,
 for a suitable suitable $\mu$-zero set $N$,
 $[\Lx 2 \mu, [H]_\sim]_{\g,2}$ can be identified with the RKHS over $X\setminus N$,
whose kernel is given by 
\begin{displaymath}
   k_\mu^\g(x,x') := 	\sum_{i\in I}  \lb_i^{\g} e_i(x) e_i(x') \, , \qquad \qquad x,x'\in X\setminus N\, .
\end{displaymath}
Since the eigenvalues of the corresponding integral operator are $\lb_i^{\g}\asymp i^{-2\g/\a}$, we conclude from 
\eqref{s-squared} and 
\eqref{ln-vs-en-2} that $e_n\bigl([\Lx 2 \mu, [H]_\sim]_{\g,2}\hookrightarrow \Lx 2\mu\bigr) \asymp  n^{-\g/\a}$.

Finally, 
using $\ran I_{k,\mu} = [H]_\sim$ and Theorem \ref{main1} the remaining assertions, namely the 
case $\b=\g=1$ as well as  the assertion for $I_{k,\mu}:H\to \Lx \infty \mu$ can be proven analogously.
\end{proofof}

Theorem \ref{interpol} essentially states that the property \eqref{interpol-result}
is passed down from the large spaces in the scale of spaces $[\Lx 2 \mu, [H]_\sim]_{s,2}$
to the smaller ones. Moreover, using the spaces on the right hand side of \eqref{interpol-represent}
instead of the interpolation spaces, it can easily be seen that the result is also true for $\g>1$. 
In addition, the representation \eqref{interpol-represent} suggests that the eigenfunctions may play 
a crucial role in determining whether  \eqref{interpol-result} holds.
In this respect note that \cite[Lemma 5.1]{MeNe10a} essentially showed the continuous embedding 
$[\Lx 2 \mu,  [H]_\sim]_{\a/2,1}\hookrightarrow \Lx \infty\mu$  provided that 
\eqref{entrop-assump} holds and that
the 
eigenfunctions are not only bounded but uniformly bounded. From this it is easy to conclude that 
$[\Lx 2 \mu,  [H]_\sim]_{\b,2}\hookrightarrow \Lx \infty\mu$ holds for all $\b\in (\a/2,1)$.
In addition, the case 
$[\Lx 2 \mu,  [H]_\sim]_{\b,2}\hookrightarrow \Lx \infty\mu$ for $\b\in (0,\a/2]$ 
can always be excluded, since 
\cite[Theorem 5.3]{StSc12a} shows that such an inclusion would imply $\sum_{i\geq 1}\lb_i^{\b}<\infty$ 
for the eigenvalues of the integral operator $T_k$, and this summability clearly contradicts 
\eqref{entrop-assump} by \eqref{s-squared} and \eqref{ln-vs-en-2}.
Summarizing, we think that understanding when 
the asymptotic equivalence \eqref{interpol-result} holds  for some $\g$ close to $\a/2$
is an interesting question for future research.

%% file: examples.tex
\section{An Example: Sobolev Spaces}\label{sec:examples}

The goal of this section is to illustrate the consequences of Lemma \ref{new-equiv} 
and Theorem \ref{main1}
by 
applying them to embeddings of the form $\id:H\to L_p(\mu)$, where $H$ is a Sobolev space, 
$X\subset \Rd$ is a suitable subset,
$\mu$ is the Lebesgue measure on $X$, and
$p\in [2,\infty]$.

We begin by recalling some basics on Sobolev spaces. 
To this end let $X\subset \Rd$ be a non-empty, open, and conncected
 subset satisfying the strong local Lipschitz condition
in the sense of \cite[p.~83]{AdFo03}. 
For $m\geq 1$ being an integer, 
we denote the classical Sobolev space on $X$ that  is defined by weak derivatives, see e.g.~\cite[p.~59-60]{AdFo03},
by $W^m(X):= W^{m,2}(X)$.

For $m>d/2$, it is well-known that the embedding $\id:W^m(X)\to C_B(X)$
is compact, see e.g.~\cite[Theorem 6.3]{AdFo03} in combination with \cite[p.~84]{AdFo03}.
Therefore, the embeddings
$\id:W^m(X)\to L_\infty(X)$ are  compact, and if $X$ has finite Lebesgue measure,
we also obtain the compactness of the 
 the embeddings
$\id:W^m(X)\to L_p(X)$, where we followed the standard notation $\Lx p X := \Lx p \mu$.
 Note that an immediate consequence of this is that the
 approximation and entropy numbers of these embeddings converge to zero.
Let us recall some results from \cite{EdTr96} that describe the asymptotic behavior of these numbers. 
To this end,
note that a consequence of Stein's extension theorem, see \cite[Theorem 5.24]{AdFo03} is that 
\begin{equation}\label{sobol-restric}
\snorm f := \inf\bigl\{  \snorm g_{W^m(\Rd)} : g\in W^m(\Rd) \mbox{ with } g_{|X} = f \bigr\}\, ,
\end{equation}
where $f\in W^m(X)$,
defines an equivalent norm on $W^m(X)$. Moreover, if for $s\in [0,\infty)$ and $p,q,\in (0,\infty]$
we write $B^s_{p,q}(\Rd)$ and $F^s_{p,q}(\Rd)$ for the  Besov  and Triebel-Lizorkin spaces
in the sense of \cite[p.~24f]{EdTr96}, then we have $B^m_{2,2}(\Rd) = F^m_{2,2}(\Rd) = W^m(\Rd)$
by \cite[p.~44 and p.~25]{EdTr96}.
By \eqref{sobol-restric} we conclude that the spaces $B^m_{2,2}(X)$ defined by restrictions as in 
\cite[p.~57]{EdTr96} satisfy
\begin{equation}\label{sobol-eq-besov}
   B^m_{2,2}(X) = W^m(X)
\end{equation}
up to equivalent norms. Moreover, \cite[p.~25]{EdTr96} shows that $F^0_{p,2}(\Rd) = L_p(\Rd)$, and by 
\cite[p.~44]{EdTr96} we find continuous embeddings $B_{p,2}^0(\Rd)  \hookrightarrow L_p(\Rd) \hookrightarrow B_{p,p}^0(\Rd)$
for all $p\in [2,\infty)$. By \eqref{sobol-restric} we conclude that 
\begin{equation}\label{besov-vs-lp}
   B_{p,2}^0(X)  \hookrightarrow L_p(X) \hookrightarrow B_{p,p}^0(X)\, .
\end{equation}
Similarly, recall that we have   continuous embeddings 
$B_{\infty,1}^0(\Rd)  \hookrightarrow L_\infty(\Rd) \hookrightarrow B_{\infty,\infty}^0(\Rd)$, see 
e.g.~\cite[p.~44]{EdTr96}, and thus we also have the continuous embeddings 
\begin{equation}\label{besov-vs-linfty}
   B_{\infty,1}^0(X)  \hookrightarrow L_\infty(X) \hookrightarrow B_{\infty,\infty}^0(X)\, .
\end{equation}
Let us now assume   that $X$ is open, connected, and bounded, and that it has a $C^\infty$-boundary. 
Moreover, we fix some  $s_1,s_2\in [0,\infty)$ and   $p_1,p_2,q_1,q_2\in (0,\infty]$
such that $s_1-s_2 > d \bigl(\frac 1 {p_1} - \frac 1 {p_2}  \bigr)_+$.
Then \cite[Theorem 2 on p.~118]{EdTr96} shows that 
\begin{equation}\label{besov-entropy}
   e_n\bigl(\id: B^{s_1}_{p_1,q_1}(X)   \to  B^{s_2}_{p_2,q_2}(X) \bigr)   \asymp n^{-(s_1-s_2)/d}\, ,
\end{equation}
and if we additionally assume that $2\leq p_1\leq p_2\leq \infty$, then \cite[p.~119]{EdTr96}
shows that 
\begin{equation}\label{besov-approx}
   a_n\bigl(\id: B^{s_1}_{p_1,q_1}(X)   \to  B^{s_2}_{p_2,q_2}(X) \bigr)   \asymp n^{-(s_1-s_2)/d + 1/p_1 - 1/p_2}\, .
\end{equation}
In particular, for $s_1=s$, $p_1=q_1=2$,  $s_2=0$, $p_2=p\in [2,\infty]$, and $q_2=q\in [1,\infty]$ 
with $s> d(\frac 12-\frac 1p)$
we obtain
\begin{align*} 
    e_n\bigl(\id: B^{s}_{2,2}(X)   \to  B^{0}_{p,q}(X) \bigr)  & \asymp n^{-s/d}\\  
  a_n\bigl(\id: B^{s}_{2,2}(X)   \to  B^{0}_{p,q}(X) \bigr)  & \asymp n^{-s/d+ 1/2-1/p}\, .
\end{align*}
By \eqref{sobol-eq-besov}, \eqref{besov-vs-lp}, and \eqref{besov-vs-linfty} we conclude that 
\begin{align}\label{entr-final}
    e_n\bigl(\id:W^m(X)   \to L_p(X) \bigr)  & \asymp n^{-m/d}\\   \label{apprx-final}
  a_n\bigl(\id: W^m(X)  \to  L_p(X) \bigr)  & \asymp n^{-m/d+ 1/2 - 1 /p} 
\end{align}
for all $m\in \N$ with $m> d(\frac 12-\frac 1p)$. In other words, the gap between the entropy and approximation
numbers is of the order $n^{1/2 - 1 /p}$. Note that for the Hilbert space case, i.e.~$p=2$, the gap 
vanishes as already observed in Section \ref{sec:carls}, while in the other extreme $p=\infty$,
the gap is of the order $n^{1/2}$.
Finally, we see by \eqref{sobol-restric}  that 
these asymptotics still hold, if we only assume that $X$ is an open, connected, and bounded subset of $\Rd$ satisfying 
the strong local Lipschitz condition.

To illustrate these findings, we now consider the linear interpolation $n$-width
mentioned in the introduction.
%
To this end, we fix an $m\in \N$ with $m>d/2$ and let 
$H=W^m(X)$ 
with equivalent norms. Then  \eqref{apprx-final} shows
\begin{displaymath}
  n^{-m/d+ 1/2-1/p} \asymp    a_{n+1}\bigl((\id :H\to L_p(X)\bigr) \leq I_n\bigl(H, L_p(X)\bigr) 
\end{displaymath}
for all $p\in [2,\infty]$.
 Here we note that in the case $p=\infty$, the lower bound 
$n^{-m/d+ 1/2} \prec I_n\bigl(H, L_\infty(X)\bigr)$ 
already follows from \eqref{in-lower-bound}.
%
Moreover,  \eqref{entr-final} in combination with Theorem \ref{main1} 
yields
\begin{equation}\label{dn-estimate}
 d_{n}\bigl(\id :H\to L_p(X) \bigr)  \asymp  n^{-m/d}
\end{equation}
for all $p\in [2,\infty]$. 
In other words, the gap of $1/2-1/p$ actually occurs between the non-linear approximation described by $d_n$
and the linear approximation described by $a_n$. Moreover, the gap is maximal for $p=\infty$ and vanishes 
in the other extreme case $p=2$.

We like to mention that \eqref{dn-estimate} appears to be new, since 
it is not contained in the list of known asymptotics compiled in \cite{Vybiral08a}. In addition, the gap
between $d_n$ and $a_n$ is solely derived from the same gap between 
$e_n$ and $a_n$, that is, from \eqref{entr-final} and \eqref{apprx-final}.
In other words, we will observe a gap between $d_n$ and $a_n$ if and only if there is a gap between
$e_n$ and $a_n$. Fortunately, the latter two quantities have been considered for various other spaces $H$ and measures $\mu$,
so that it should be possible to compile a list of cases, in which the gap occurs.

For convenience, the following corollary summarizes our findings for sufficiently large subspaces of $W^m(X)$.
It particularly applies to kernels of many 
standard Gaussian processes, such as the (iterated) Brownian motion and -bridge,
see e.g.~the numerical example in \cite{SaSc16a}.

\begin{corollary}\label{sobol-cor}
Let $X\subset \Rd$ be an open, connected, and bounded subset satisfying the strong Lipschitz condition.
Moreover, let $H$ be an RKHS over $X$ with kernel $k$
such that $H\hookrightarrow W^m(X)$ for some integer $m>d/2$.
Assume, in addition, that 
\begin{displaymath}
   e_{n}\bigl(I_{k,\mu} :H\to L_2(X) \bigr) \asymp  n^{-m/d}
\end{displaymath}
holds. Then we have 
\begin{displaymath}
   d_n\bigl(H, \Lx \infty \mu\bigr)\asymp n^{-m/d} \qquad \qquad  \mbox{ and } 
\qquad \qquad n^{-m/d + 1/2} \prec  I_n\bigl(H, L_\infty(\mu)\bigr) \, .
\end{displaymath}
In addition, if $H=W^m(X)$ with equivalent norms, then, for all $p\in [2,\infty]$, we have 
\begin{displaymath}
   d_n\bigl(H, \Lx p \mu\bigr)\asymp n^{-m/d} \qquad \qquad  \mbox{ and } 
\qquad \qquad n^{-m/d + 1/2 - 1/p} \asymp  a_n\bigl(H, L_p(\mu)\bigr)\, .
\end{displaymath}
\end{corollary}

\begin{proofof}{Corollary \ref{sobol-cor}}
 We first note that the sequence of estimates 
\begin{displaymath}
 n^{-m/d} 
\asymp   
e_{n}\bigl(I_{k,\mu} :H\to L_2(X) \bigr)
\prec
e_{n}\bigl(I_{k,\mu} :H\to L_\infty(X) \bigr)
\prec 
e_{n}\bigl(I_{k,\mu} :W^m\to L_\infty(X) \bigr)
\asymp 
 n^{-m/d} 
\end{displaymath}
yields $e_{n}\bigl(I_{k,\mu} :H\to L_\infty(X) \bigr) \asymp n^{-m/d}$, and therefore Corollary \ref{gap-1}
shows the first two assertions. 
The second set of asymptotic equivalences immediately follows from our findings 
in the text above together with the muliplicativity of the approximation numbers.
\end{proofof}

%